\documentclass[a4paper,12pt]{article}
\usepackage{hyperref}
\usepackage[english]{babel}
\usepackage[T1]{fontenc}
\usepackage[utf8]{inputenc}
\usepackage{amsfonts}
\usepackage{amssymb}
\usepackage{amsmath}
\usepackage{graphicx}
\usepackage{cite}
\usepackage{epsfig}
\usepackage{psfrag}
\usepackage{tikz}
\usetikzlibrary{matrix,arrows,decorations.pathmorphing,positioning} 
\usepackage{url}
\usepackage{amsthm}
\usepackage{mathpartir}
\usepackage{enumerate}
\usepackage{tikz-cd}

\newtheorem{proposition}{Proposition}[section]
\newtheorem{definition}{Definition}[section] 
  
\newtheorem{corollary}{Corollary}[section] 
\newtheorem{remark}{Remark}[section]  

\date{}

\title{\bf Homotopies in Multiway (Non-Deterministic) Rewriting Systems as $n$-Fold Categories}

\author{Xerxes D. Arsiwalla$^{1, 4, }$\footnote{Corresponding Author: \url{x.d.arsiwalla@gmail.com}}  {\,\,}   Jonathan Gorard$^{2, 4, }$\footnote{\url{jg865@cam.ac.uk}}  {\,\,} Hatem Elshatlawy$^{3, 4, }$\footnote{\url{hatem.elshatlawy@rwth-aachen.de}} \\
{}  \\
{\it \small $^{1}$Pompeu Fabra University, Barcelona, Spain}\\ 
{\it \small $^{2}$University of Cambridge, Cambridge, United Kingdom}\\
{\it \small $^{3}$RWTH Aachen University, Aachen, Germany}\\
{\it \small $^{4}$Wolfram Research, USA}    
}

\begin{document}
\maketitle

\begin{abstract}
We investigate algebraic and compositional properties of abstract multiway rewriting systems, which are archetypical structures underlying the formalism of the \textit{Wolfram model}. We demonstrate the existence of higher homotopies in this class of rewriting systems, where homotopical maps are induced by the inclusion of appropriate rewriting rules taken from an abstract \textit{rulial space} of all possible such rules. Furthermore, we show that a multiway rewriting system with homotopies up to order $n$ may naturally be formalized as an $n$-fold category, such that (upon inclusion of appropriate inverse morphisms via invertible rewriting relations) the infinite limit of this structure yields an ${\infty}$-groupoid. Via Grothendieck's homotopy hypothesis, this ${\infty}$-groupoid thus inherits the structure of a formal homotopy space. We conclude with some comments on how this computational framework of homotopical multiway systems may potentially be used for making formal connections to homotopy spaces upon which models relevant to  physics may be instantiated.
\end{abstract}

\clearpage

\tableofcontents

%\clearpage

\section{Introduction}

The \textit{Wolfram model}\cite{Wolfram2020} \cite{Wolfram2002a}  (see also \cite{Gorard2020} \cite{Gorard2020a} \cite{Gorard2020c}  \cite{arsiwalla2021pregeometric}) is an explicitly computational framework based upon abstract rewriting that attempts to formalize the fundamental mathematical structures underlying space, time and matter. The model seeks to capture the ways in which simple abstract rewriting rules may be composed, so as to yield more complex structures (suggestively titled ``universes'') that admit certain emergent ``laws of physics''. This general approach of seeking a constructivist formalization of certain fundamental structures has strong parallels in recent work in the foundations of both physics and mathematics, for instance in relation to synthetic geometry and cohesive homotopy type theory\cite{Schreiber2012}\cite{Schreiber2013a}\cite{Shulman2016}\cite{Ahrens2021}. The general framework of the Wolfram model eschews the rigid continuum description of spacetime in terms of Lorentzian manifolds in favor of a more rudimentary description of the combinatorics and topology of the intrinsic causal relationships between events, which, in turn, makes manifest some of the computational architecture that underlies contemporary physics. The archetypical structures that appear within this framework are so called \textit{multiway systems} - non-deterministic abstract rewriting systems equipped with a notion of causal structure. The ``multiway'' moniker here designates the fact that all permissible applications of the rewriting rules are instantiated in all possible orderings, leading to multiple chains of rewrite terms that are partially-ordered by causality. Depending upon the precise interpretation of these terms (or \textit{states}), Wolfram model multiway systems may be realized as rewriting systems over hypergraphs, character strings, ZX-diagrams\cite{Gorard2020c}\cite{Gorard2021a} (as formalized by Coecke and Duncan\cite{Coecke2009a}), or  string diagrams\cite{Gorard2021b}.

In this paper, we consider \textit{abstract} multiway rewriting systems, i.e. systems in which the terms/states are purely abstract entities, since doing so will be more convenient for studying the algebraic and compositional properties of such systems at the highest possible level of generality; these general results may then be easily specialized or extrapolated to particular classes of rewriting systems. Specifically, our analysis reveals a rich homotopical structure admissible on these abstract multiway rewriting systems, where these homotopies are induced by the inclusion of certain ``higher-order'' rules taken from an abstract \textit{rulial space}, designating the space of all possible rewriting rules of a given signature (formally, at least for the case of hypergraphs, a monoidal category of cospans). Although a generic abstract rewriting system may be thought of as simply being an $F$-coalgebra for the power set functor  ${\mathcal{P}}$  (which, at least in the case of hypergraph rewriting, is further equipped with a natural symmetric monoidal structure), we show here that a multiway rewriting system equipped with homotopies up to order $n$ may elegantly be formalized as an $n$-fold category. Furthermore, upon including inverse morphisms (via the inclusion of invertible rewriting relations), the infinite limit of this structure naturally yields an ${\infty}$-groupoid, which therefore inherits the structure of a formal homotopy space via Grothendieck's homotopy hypothesis \cite{Arsiwalla2020}.  For specific models of ${\infty}$-categories, these structures correspond precisely to topological spaces, hence allowing one to make formal connections between the abstract computational framework of Wolfram model multiway systems and the topological spaces (such as Lorentzian manifolds, or projective Hilbert spaces) upon which laws of physics are traditionally instantiated.

\section{Preliminaries of Wolfram Model Multiway Systems}
\label{sec:Section1}

This section provides an overview of the basic abstract rewriting constructions that constitute the Wolfram model, with a particular emphasis on their non-deterministic aspects, as captured via \textit{multiway rewriting systems}. We present both combinatorial and compositional descriptions of such systems, as well as an abstract description of their \textit{rulial space} (i.e. the space of all possible rewriting rules) for the case of (hyper)graph rewriting, by means of a double-pushout (DPO) rewriting formalism over (selective) adhesive categories. The contents of this section serve merely as a preliminary introduction to the known constructions of Wolfram model multiway rewriting systems; the informed reader is encouraged to proceed directly to Section \ref{sec:Section2}, in which the main arguments of this paper are presented.

\subsection{Combinatorial \& Compositional Descriptions of Multiway  Systems}

Combinatorially, a multiway rewriting system is simply a directed, acyclic graph of states, determined by abstract rewriting rules that inductively generate a (potentially infinite) \textit{multiway evolution graph}, together with a partial order on the rewrite rule applications, determined by their causal structure. This is straightforward to formalize in terms of \textit{abstract rewriting systems} \cite{Bezem2003} in which the underlying rewrite relation is not (necessarily) confluent\cite{Huet1980}.

\begin{definition}
An ``abstract rewriting system'' (or ``ARS'') is a set, denoted $A$ (whose elements are known as ``objects''), equipped with a binary relation, denoted ${\to}$ (known as the ``rewrite relation'').
\end{definition}

\begin{figure}[ht]
\centering
\includegraphics[width=0.695\textwidth]{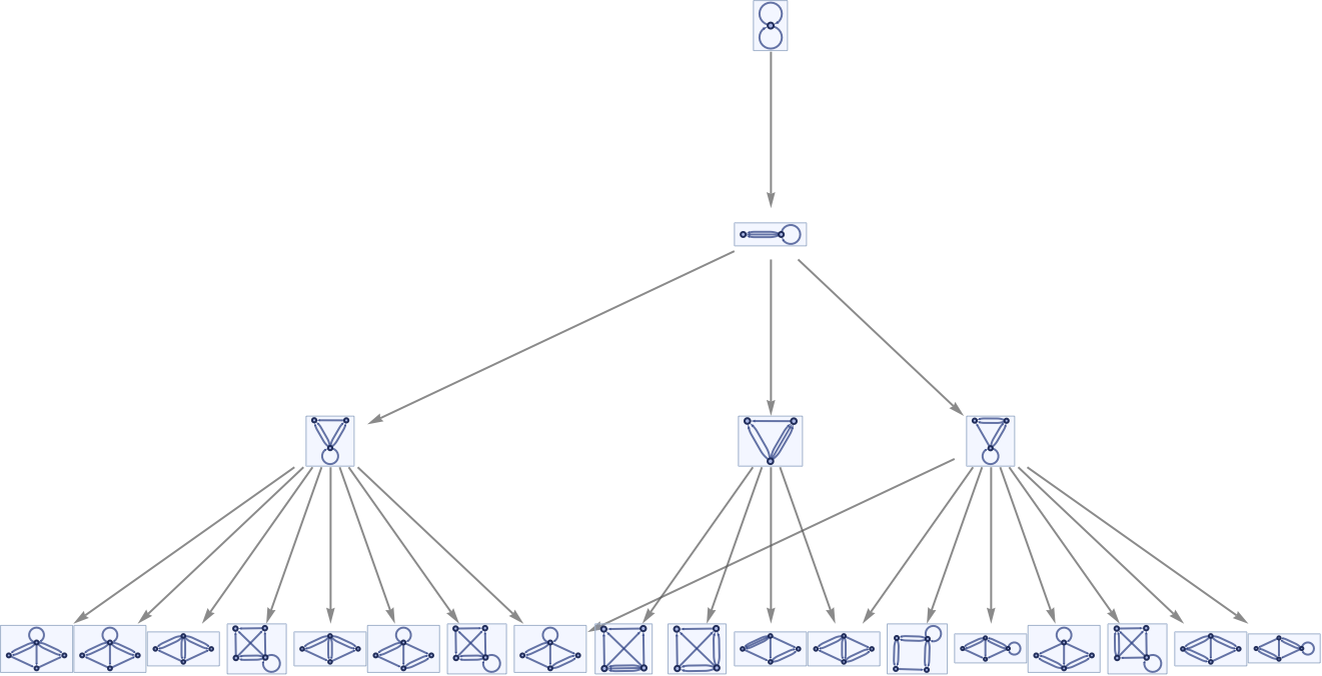}
\caption{The multiway evolution graph corresponding to the first 3 steps in the non-deterministic evolution history of the hypergraph substitution rule ${\left\lbrace \left\lbrace x, y \right\rbrace, \left\lbrace y, z \right\rbrace \right\rbrace \to \left\lbrace \left\lbrace w, y \right\rbrace, \left\lbrace y, z \right\rbrace, \left\lbrace z, w \right\rbrace, \left\lbrace x, w \right\rbrace \right\rbrace}$, starting from a simple double self-loop initial condition ${\left\lbrace \left\lbrace 0, 0 \right\rbrace, \left\lbrace 0, 0 \right\rbrace \right\rbrace}$.}
\label{fig:Figure1}
\end{figure}

\begin{definition}
A ``multiway evolution graph'', denoted ${G_{multiway}}$, is a directed, acyclic graph corresponding to the evolution of a (generically non-confluent) abstract rewriting system ${\left( A, \to \right)}$, in which the set of vertices corresponds to the set of objects ${V \left( G_{multiway} \right)}$, and in which the directed edge ${A \to B}$ exists in ${E \left( G_{multiway} \right)}$ if and only if there exists an application of the rewrite relation that transforms object $A$ to object $B$.
\end{definition}
Hence, directed edges will connect vertices $A$ and $B$ in ${G_{multiway}}$ if and only if ${A \to B}$ in the underlying rewriting system, and a directed path will connect vertices $A$ and $B$ if and only if ${A \to^{*} B}$, where ${\to^{*}}$ denotes the reflexive transitive closure of ${\to}$, i.e if and only if there exists a finite rewrite sequence of the form:

\begin{equation}
A \to A^{\prime} \to A^{\prime \prime} \to \cdots \to B^{\prime} \to B.
\end{equation}
An example of a multiway evolution graph for a simple Wolfram model evolution (i.e. a non-deterministic hypergraph rewriting system) is shown in Figure \ref{fig:Figure1}.

In the most general case, in which the rewrite relation ${\to}$ is treated as an indexed union of sub-relations, i.e. ${\to = \to_1 \cup \to_2 \cup \dots}$, with label set ${\Lambda}$ for the indices, and in which the set of objects $A$ is arbitrary (i.e. its elements could represent graphs, hypergraphs, string diagrams, terms, character strings, etc.), the resulting labeled abstract rewriting system ${\left( A, \Lambda, \to \right)}$ permits an elegant compositional description in terms of $F$-coalgebras. Specifically, note that the system ${\left( A, \Lambda, \to \right)}$ is now simply a bijective function from $A$ to a subset of the power set of $A$ indexed by ${\Lambda}$, i.e. ${\mathcal{P} \left( \Lambda \times A \right)}$:

\begin{equation}
p \mapsto \left\lbrace \left( \alpha, q \right) \in \Lambda \times A : p \to^{\alpha} q \right\rbrace.
\end{equation}
Recall now that an $F$-coalgebra for an endofunctor ${F : \mathbf{C} \to \mathbf{C}}$ consists of an object $A$ in ${\mathrm{ob} \left( \mathbf{C} \right)}$ equipped with a morphism ${\alpha : A \to F A}$ in ${\mathrm{hom} \left( \mathbf{C} \right)}$, hence denoted ${\left( A, \alpha \right)}$. Thus, since the power set construction on ${\mathbf{Set}}$ is a covariant endofunctor ${\mathcal{P} : \mathbf{Set} \to \mathbf{Set}}$, we see that the abstract rewriting system ${\left( A, \to \right)}$ consists of an object $A$ equipped with an additional morphism of ${\mathbf{Set}}$, namely the rewrite relation ${\to}$:

\begin{equation}
\to : A \to \mathcal{P} A.
\end{equation}

The causal structure of the abstract rewriting system also endows the multiway system with an additional partial order relation on the rewrite applications themselves\cite{Gorard2020}, which can in turn be described compositionally in terms of a partial monoidal structure\cite{Gorard2020c}\cite{Gorard2021a}\cite{Gorard2021b}. However, for the purposes of the present paper, the causal structure can (for the most part) be safely ignored.

\subsection{Rulial Space and Double-Pushout (DPO) Rewriting}

In a large class of cases, such as (hyper)graph rewriting, open graph rewriting, string diagram rewriting (including the ZX-calculus), etc., the rewrite relation ${\to}$ may be further specified as a double-pushout (DPO) rewriting system; for instance, this was formalized by Dixon and Kissinger\cite{Dixon2013} for the case of open graphs using a selective adhesive category of open spans, and the same formalism was subsequently applied to both ZX-diagrams and arbitrary hypergraph rewriting systems as special cases\cite{Gorard2020c}\cite{Gorard2021a}. Specifically, the rewriting rules are defined as spans of monomorphisms ${\rho}$ of the form:

\begin{equation}
\rho = \left( l : K \to L, r : K \to R \right),
\end{equation}
with the left- and right-hand-sides of the rule specified by objects $L$ and $R$, respectively, and with object $K$ designating the \textit{interface graph}. A \textit{match} for the rule ${\rho}$ within an object $G$ is simply a morphism ${m : L \to G}$, with the rule being \textit{applicable} if and only if there exists a pair of pushout diagrams of the form \cite{Habel2001}:

\begin{equation}
\begin{tikzcd}
L \arrow[d, "m"] & K \arrow[l, "l"] \arrow[d, "k"] \arrow[r, "r"] & R \arrow[d, "n"]\\
G & D \arrow[l, "f"] \arrow[r, "g"] & H
\end{tikzcd}.
\end{equation}
These rules are defined and applied in the context of an \textit{adhesive category}\cite{Lack2004}, in which every pushout along a monomorphism satisfies the \textit{van-Kampen square condition}, such that, for every commutative diagram of the form:

\begin{equation}
\begin{tikzcd}
B^{\prime} \arrow[ddd, "f_{h}^{\prime}"] \arrow[dr, "h_B"] & & & A^{\prime} \arrow[lll, "g_h"] \arrow[dl, "h_A"] \arrow[ddd, "f_h"]\\
& B \arrow[d, "f^{\prime}"] & A \arrow[l, "g"] \arrow[d, "f"] &\\
& D & C \arrow[l, "g^{\prime}"] &\\
D^{\prime} \arrow[ur, "h_D"] & & & C^{\prime} \arrow[lll, "g_{h}^{\prime}"] \arrow[ul, "h_C"]
\end{tikzcd},
\end{equation}
such the following subdiagrams are both pullbacks:

\begin{equation}
\begin{tikzcd}
B^{\prime} \arrow[d, "h_B"] & A^{\prime} \arrow[l, "g_h"] \arrow[d, "h_A"]\\
B & A \arrow[l ,"g"]
\end{tikzcd}, \qquad \text{ and } \qquad
\begin{tikzcd}
A \arrow[d, "f"] & A^{\prime} \arrow[l, "h_A"] \arrow[d, "f_h"]\\
C & C^{\prime} \arrow[l, "h_C"]
\end{tikzcd},
\end{equation}
the pushouts and pullbacks satisfy a compatibility condition. In fact, \textit{selective adhesive} categories\cite{Dixon2013} allow DPO rewriting to be defined in an even broader class of cases, in which one is working within a full subcategory ${\mathbf{C}^{\prime}}$ of an adhesive category ${\mathbf{C}}$, with embedding functor ${S : \mathbf{C}^{\prime} \to \mathbf{C}}$, with the only condition being that $S$ preserves monomorphisms.

This construction allows us to formalize the notion of \textit{rulial space}, i.e. the space of all possible rewriting rules of a given signature, at least for the case of Wolfram model hypergraph rewriting systems (and systems such as string diagram rewriting systems that are known to be special cases):

\begin{definition}
The ``rulial space'' of Wolfram model systems is the category of cospans of a DPO rewriting system, defined over the selective adhesive category of Wolfram model hypergraphs.
\end{definition}
The rulial space of Wolfram model systems functorially acquires the structure of a selective adhesive category. Indeed, as a consequence of the concurrency and parallelism theorems of algebraic graph transformation theory\cite{Ehrig2006}, the rulial space also inherits a natural monoidal structure\cite{Gorard2020c} (which is, in turn, inherited by all Wolfram model multiway systems, since they are obtained by the composition of certain rules in rulial space).

Even though the rewriting formulation above borrows from DPO rewriting as per  \cite{Dixon2013}, it is worth noting other recent developments in rewriting theory, that  will no doubt be useful for our model in future work. We point the reader to the following relevant works: (i) mathematical foundations of  string diagram rewriting, reported in  \cite{bonchi2020stringI},  \cite{bonchi2021stringII}; (ii) higher dimensional rewriting theory \cite{Burroni1993},  \cite{Guiraud2019}; (iii) algebraic graph transformations and compositional systems  \cite{heckel2006graph},  \cite{behr2021compositionality}.

\section{Homotopical Multiway Systems as $n$-Fold Categories}
\label{sec:Section2}

Given the compositional setup outlined above, we now address explicit constructions  of higher homotopies upon abstract multiway systems, and introduce a systematic algorithm for identifying rewriting rules resulting in these homotopies. Subsequently, we prove that a multiway system equipped with homotopies up to order $n$ may be formalized as an $n$-fold category, such that the infinite limit of such a higher-order multiway system yields an ${\infty}$-groupoid (upon the admission of invertible rewriting rules).

\subsection{A Computational Framework for Constructing Homotopical Multiway  Systems}

\begin{figure}[ht]
\centering
\includegraphics[width=0.495\textwidth]{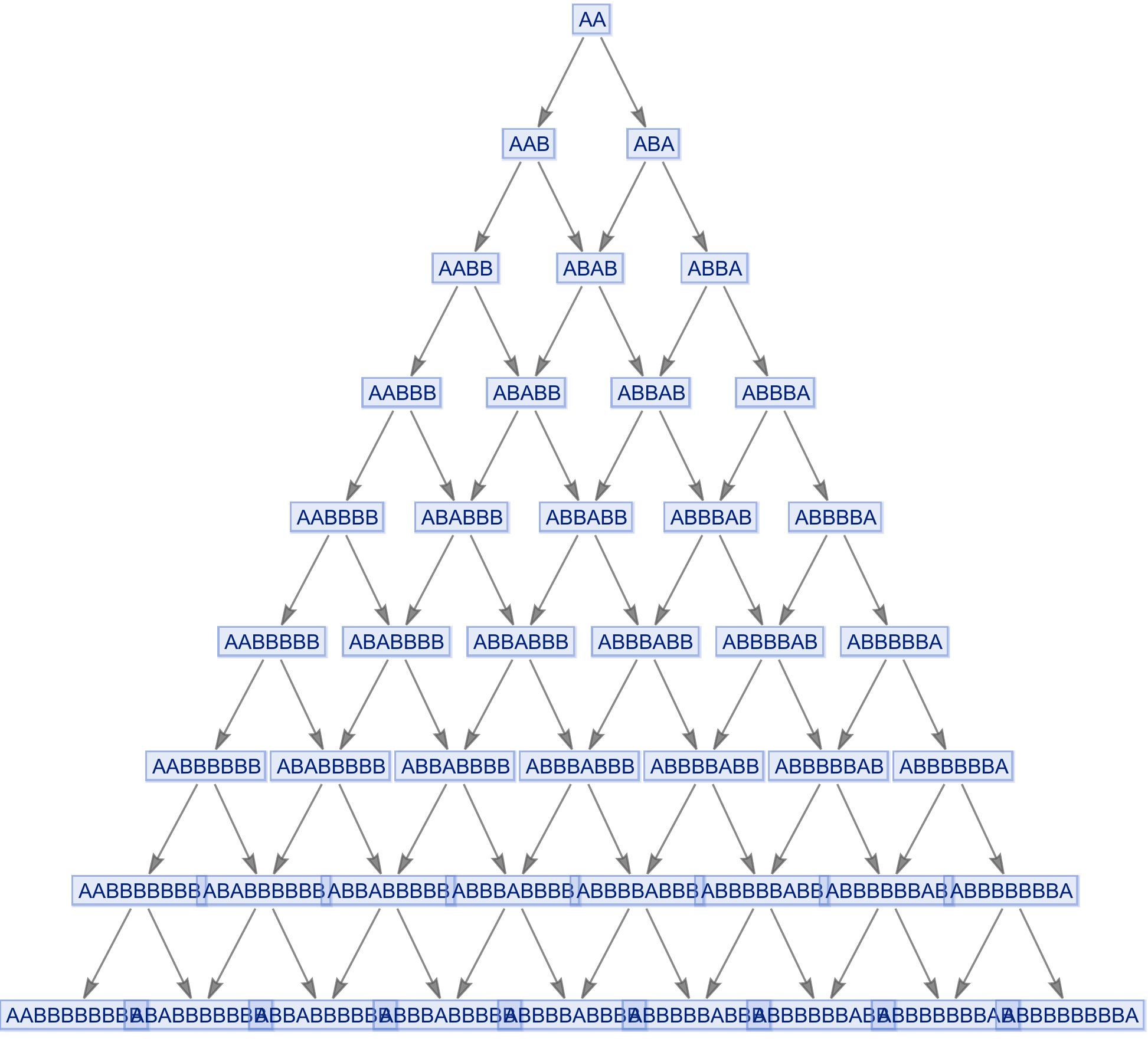}
\caption{The multiway evolution graph corresponding to the first 8 steps in the non-deterministic evolution history of the string substitution rule ${A \to AB}$, starting from the two-character initial condition $AA$.}
\label{fig:Figure2}
\end{figure}

\begin{figure}[ht]
\centering
\includegraphics[width=0.3\textwidth]{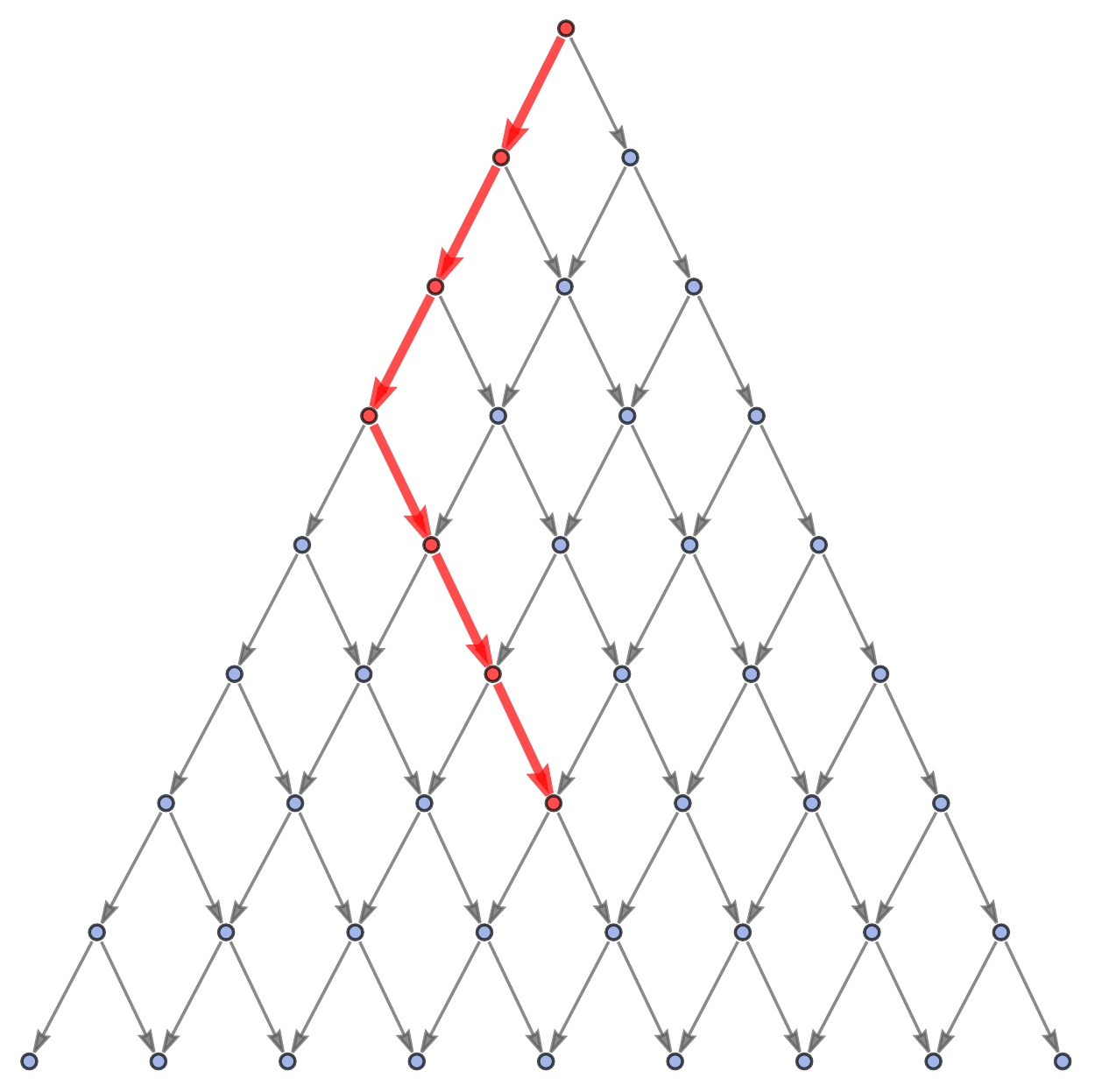}\hspace{0.1\textwidth}
\includegraphics[width=0.3\textwidth]{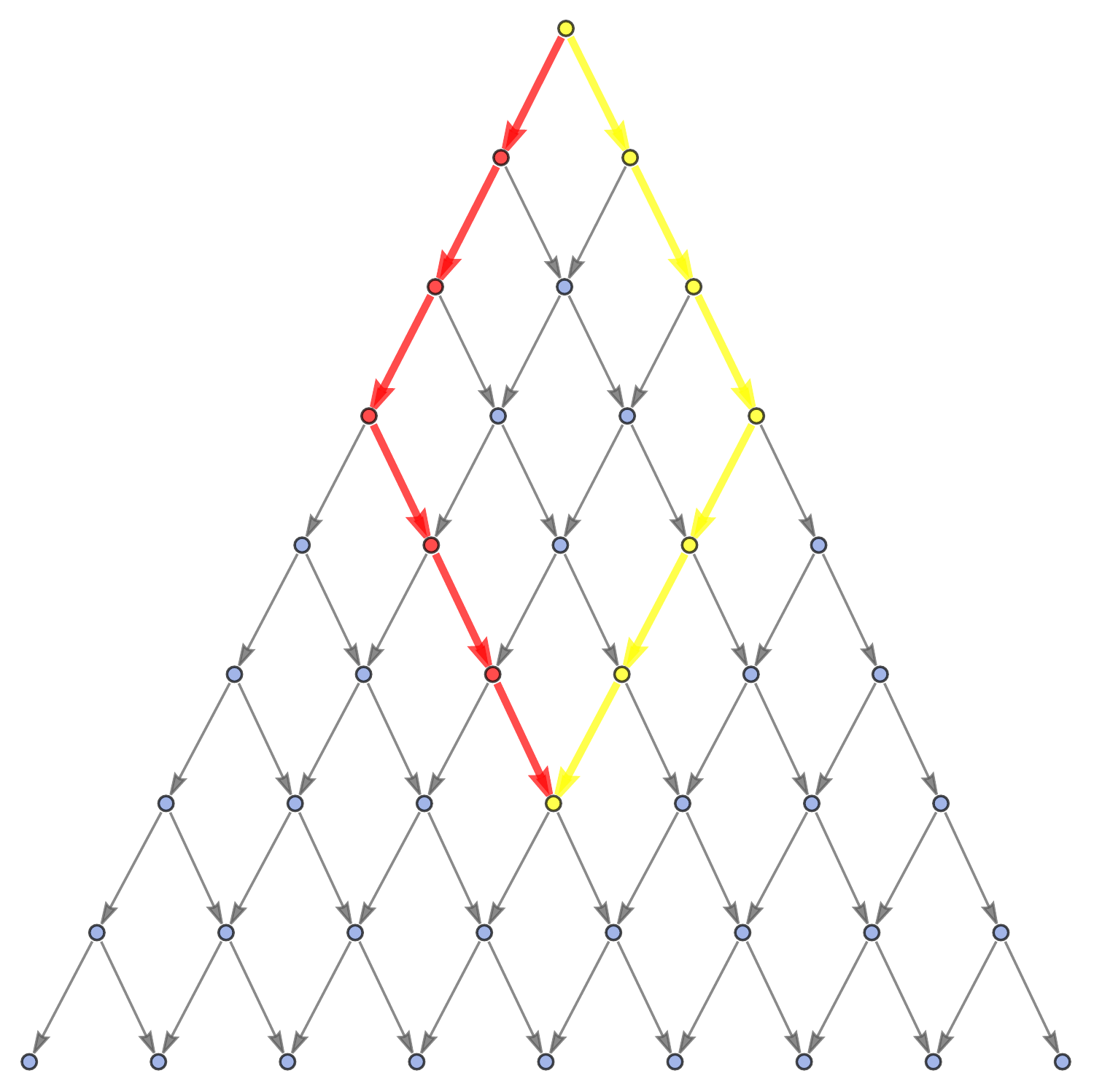}
\caption{On the left, a multiway evolution graph with a red highlighted path between vertices ${AA}$ and ${ABBBABBB}$, corresponding to a proof of the proposition ${AA \to ABBBABBB}$, subject to the string substitution axiom ${A \to AB}$. On the right, the same multiway evolution graph showing multiple highlighted paths (red and yellow) between vertices ${AA}$ and ${ABBBABBB}$, illustrating the existence of multiple proofs of the proposition ${AA \to ABBBABBB}$.}
\label{fig:Figure3}
\end{figure}

To avoid us needing to deal with the added complexities of full hypergraph rewriting systems, let us consider instead a very simple rewriting system on character strings, as specified by the rule ${A \to AB}$, generating the multiway evolution graph shown in Figure \ref{fig:Figure2}. We have deliberately selected a very minimal example, although the basic algorithm described here applies in principle to \textit{all} classes of multiway systems, as discussed below. Every path in this multiway evolution graph corresponds to a proof of a proposition, as shown on the left of Figure \ref{fig:Figure3} for the particular case of a proof of the proposition that ${AA \to ABBBABBB}$, subject to the axiom ${A \to AB}$. A key feature of multiway systems is that one can consider multiple paths connecting the same pair of vertices, which may be interpreted type-theoretically as corresponding to the existence of multiple proofs of the same proposition, as shown on the right of Figure \ref{fig:Figure3}.

\begin{figure}[ht]
\centering
\includegraphics[width=0.9\textwidth]{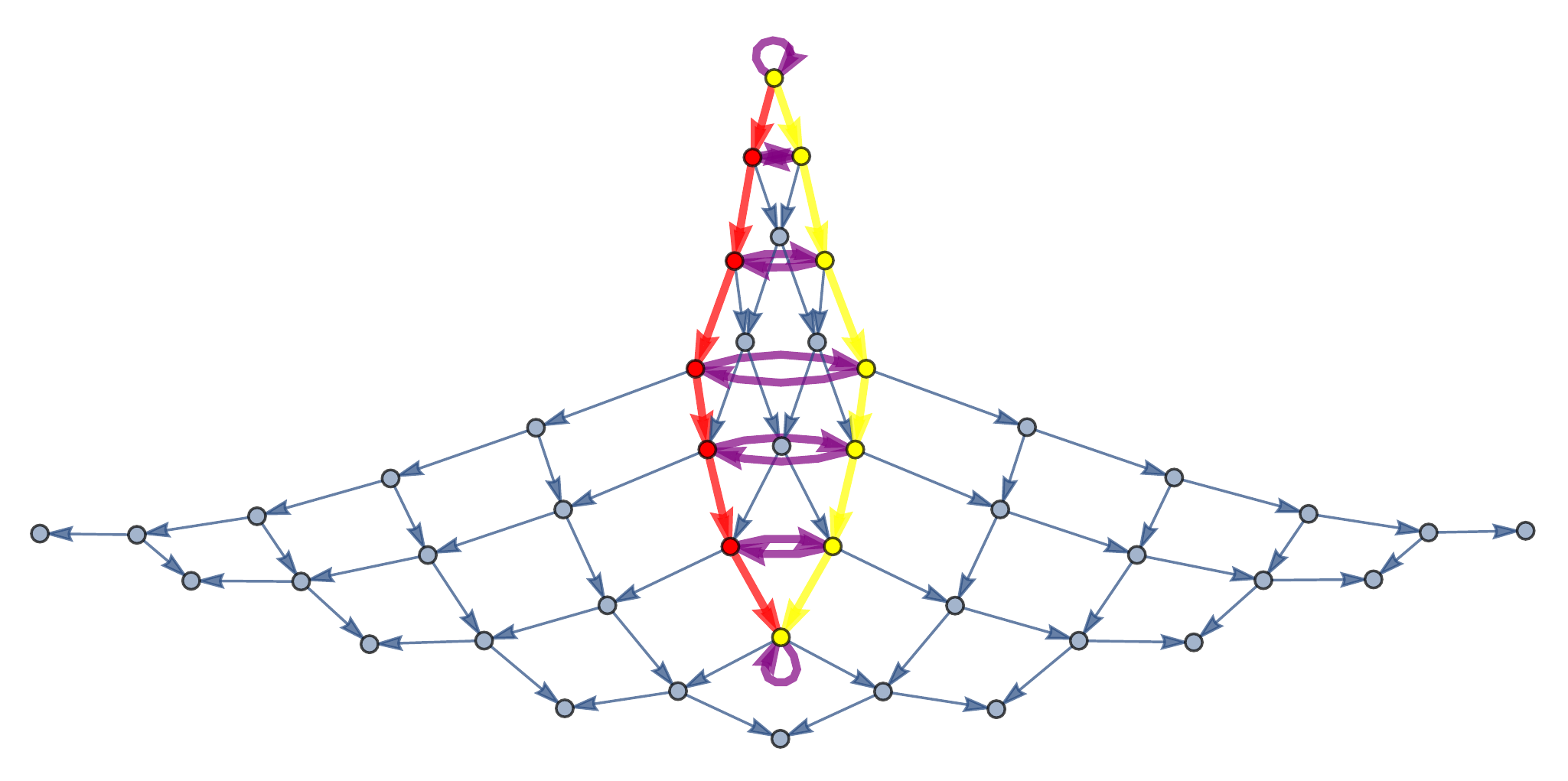}
\caption{A multiway evolution graph with purple highlighted paths between vertices on one (red) path and corresponding vertices on the other (yellow) path, and back again, interpreted as homotopic maps between the two associated proofs of the proposition ${AA \to ABBBABBB}$.}
\label{fig:Figure4}
\end{figure}

\begin{figure}[ht]
\centering
\includegraphics[width=0.4\textwidth]{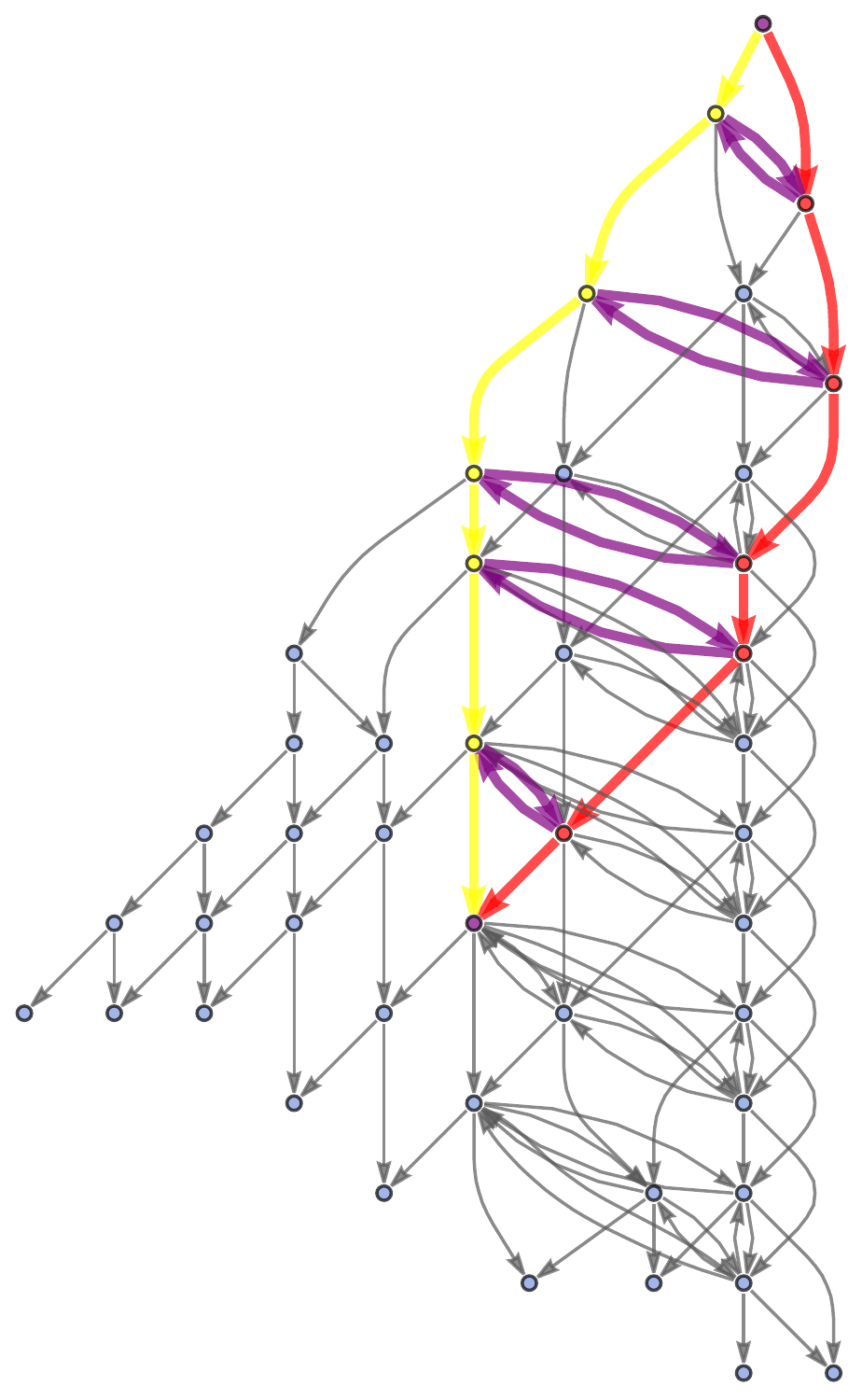}
\caption{A multiway evolution graph with explicit rules added to map vertices on one (red) path to corresponding vertices on the other (yellow) path, and back again, interpreted as homotopic maps between the two associated proofs of the proposition ${AA \to ABBBABBB}$, shown in purple.}
\label{fig:Figure5}
\end{figure}

A key feature of homotopy type theory is the ability to describe proofs of equivalences between proofs in terms of homotopies between paths\cite{Program2013}. In the discrete setup of a multiway evolution graph, a concrete realization of the homotopy map between paths can be given by simply introducing bi-directional edges mapping vertices on one path to corresponding vertices on the other path, and back again, as shown in Figure \ref{fig:Figure4}. However, this approach to introducing homotopy maps is somewhat ad hoc; it is far more natural to introduce these bi-directional relationships between multiway vertices as a new set of rewrite relations that can be appended to the original multiway system rules, as shown in Figure \ref{fig:Figure5}. For the particular system at hand, these additional rules are:

\begin{multline}
\hspace{2cm}  \left\lbrace AAB \to ABA, \,  AABB \to ABBA, \, AABBB \to ABBBA, \right.\\  
\left. ABABBB \to ABBBAB, \,  ABBABBB \to ABBBABB \right\rbrace  \hspace{3cm}
\end{multline}
and their respective inverses (which have significance for the forthcoming constructions of groupoids):

\begin{multline}
\hspace{2cm} \left\lbrace ABA \to AAB, \, ABBA \to AABB, \, ABBBA \to AABBB, \right.\\
\left. ABBBAB \to ABABBB, \, ABBBABB \to ABBABBB \right\rbrace \hspace{3cm}
\end{multline}

This general algorithmic procedure for introducing homotopy maps between paths in an arbitrary multiway system may now be iterated, so as to introduce higher homotopy maps between homotopies; for instance, order-3 homotopies (3-cells) between order-2 homotopies (2-cells) and so on (more of this will be reported in forthcoming work  \cite{arsiwalla2021pregeometric}). 

%in Figure \ref{fig:Figure6} (using the ad hoc method) and Figure \ref{fig:Figure7} (using the natural method).

%\begin{figure}[ht]
%\centering
%\includegraphics[width=0.895\textwidth]{ACTHomotopy7}
%\caption{A multiway evolution graph with orange highlighted paths corresponding to order-2 homotopies between the order-1 homotopies shown in purple between the paths shown in red and yellow. Here, the 2-cells consist of a pair of 2-squares glued together (as opposed to the more general case of a full cube with additional horizontal arrows between the squares).}
%\label{fig:Figure6}
%\end{figure}
%
%\begin{figure}[ht]
%\centering
%\includegraphics[width=0.7\textwidth]{ACTHomotopy8}
%\caption{A multiway evolution graph with explicit rules added to enact the order-2 homotopies between the order-1 homotopies between the red and yellow paths (with order-1 homotopies shown in purple and order-2 homotopies shown in orange). Here, the 2-cells consist of a pair of 2-squares glued together (as opposed to the more general case of a full cube with additional horizontal arrows between the squares).}
%\label{fig:Figure7}
%\end{figure}

\clearpage

\subsection{Higher Categorical Formulation of Homotopical Multiway Systems}

We now attempt to formalize the above constructions in the language of higher category theory. Starting with the order-2 homotopies (2-cells or 2-morphisms) between paths (1-morphisms) we make the following proposition: 
\begin{proposition}
The multiway rewriting system shown in Figure \ref{fig:Figure5} is an order-2 homotopy rewriting system with 2-morphisms between paths, thus yielding a double category.
\label{p3.1}
\end{proposition}
Let us first recall the standard definition of a double category.

\begin{definition}
A double category $D$, denoted  \begin{tikzcd}  D_1 \arrow[r, shift left]   \arrow[r, shift right]  & D_0   \end{tikzcd}  is defined in terms of the following conditions:
\begin{enumerate}[(i)]
\item
The objects of $D$ are the objects of $D_0$ 
\item
$D$ has vertical morphisms, which are the morphisms of ${D_0}$ 
\item
Additionally, $D$ also has horizontal morphisms, which are the objects of ${D_1}$  
\item
Finally, the 2-morphisms of $D$ (also referred to as squares or 2-cells) are the morphisms of ${D_1}$ 
\end{enumerate}
\end{definition}
Here, a 2-cell in $D$ can be represented by the following commutative square:
\begin{equation}
\begin{tikzcd}[column sep={2cm,between origins},row sep={2cm,between origins}]
a \dar{l} \rar["f",""{name=U, below}""] & c \dar{m} \\
b \rar["g" below,""{name=D}""] & d \\
\arrow[Rightarrow, from=U, to=D,"\phi"]
\end{tikzcd} 
\end{equation}
where $a$, $b$, $c$, $d$ are objects; $l$, $m$ are vertical arrows; $f$, $g$ are horizontal arrows; and $\phi$ denotes the 2-cell. 

Note that a double category is an internal category in {\bf Cat}. Vertical composition in a double category is given by composition in the categories $D_0$ and $D_1$, while horizontal composition is given by composition on \begin{tikzcd}  D_1 \arrow[r, shift left]   \arrow[r, shift right]  & D_0   \end{tikzcd} itself, by virtue of it being a category internal to {\bf Cat}. This set-up will be useful for proving the desired proposition.

{\it Proof of Proposition \ref{p3.1}}: 
Following the functorial description provided in Section \ref{sec:Section1}, any abstract  multiway system may be thought of as a category ${\cal M}_0$ whose  objects are given by states (rewrite terms) of the multiway system, and whose paths of causally-ordered rewrite chains refer to morphisms of this category. Composition is given by concatenation of directed paths. Transitive closure along directed paths ensures associativity and self-loops on nodes (trivially obtained via identity rewrite rules) ensure identity laws of a category.

From the definition of a double category given above, ${\cal M}_0$ is precisely the category $D_0$, with objects being the vertices of a multiway evolution graph (a prototypical example being given in Figure \ref{fig:Figure2}),  and paths formed by multiway evolution (also shown in Figure \ref{fig:Figure3}) corresponding to vertical morphisms in $D_0$. 
The horizontal arrows in Figure \ref{fig:Figure5}, colored in purple, are induced by the introduction of additional rules to the rewriting system and serve the purpose of connecting (vertical) paths in the multiway system. Hence, these can be identified as objects of another category, which we denote as ${\cal M}_1$ (${\cal M}_1$ will be identified with $D_1$ once we specify its morphisms).
Imposing commutativity over squares formed by the vertical and horizontal arrows of the multiway system gives us the morphisms of ${\cal M}_1$; these squares, exemplified in Figure \ref{fig:Figure5}, are 2-morphisms of the double category   \begin{tikzcd}  {\cal M}_1 \arrow[r, shift left]   \arrow[r, shift right]  & {\cal M}_0   \end{tikzcd}. 

Given this double category, the 2-morphism in the multiway evolution graph between the highlighted red and yellow paths in Figure \ref{fig:Figure5} is precisely given by vertical compositions of squares. These are the 2-cells of the new multiway rewriting system, enhanced with additional rules, thus making it a \textit{2-homotopy rewriting system}. (Note that the triangles at the start and end points of the 1-paths are also squares with an identity morphism at those points.)
\qed \\  \; 

This leads us directly to the following corollary: 
\begin{corollary}
\label{c3.1}
A 2-homotopy multiway system extends to a double groupoid upon admitting invertible rules to those used in the multiway construction (including the homotopy rules).
\end{corollary}

{\it Proof of Corollary \ref{c3.1}}: 
Adding invertible rules to the multiway system equips every rewrite arrow with an inverse arrow, such that their composition is the identity map on either rewrite term. Given an order-2 homotopy  multiway system, the presence of invertible arrows implies that all 1-morphisms and 2-morphisms are isomorphisms. This extends the multiway system in Proposition \ref{p3.1} to a double groupoid.  \qed 

Furthermore, as shown in the previous section, with appropriate additional rules, one can continue the above construction inductively to obtain rewriting systems with progressively higher homotopies between paths. This yields the following proposition:

\begin{proposition}
So long as additional rewrite rules between cells, up to order $n - 1$, are admissible, then the multiway system shown in Figure \ref{fig:Figure5} can be enhanced to yield an $n$-homotopy rewriting system with $n$-morphisms between paths, thus yielding an $n$-fold category.
\label{p3.2}
\end{proposition}

{\it Proof of Proposition \ref{p3.2}}: 
The proof of this proposition follows by a simple inductive construction. 

A double category itself is a  special case of an $n$-fold category. In Proposition \ref{p3.2}, we proved this proposition for the $n = 2$ case.  

For the $n = 3$ case, we need to construct an order-3 homotopy rewriting system equipped with 3-morphisms, which are cubes bounded by squares of the given double category. This can be done by adding specific rewrite rules that give arrows between paths (1-cells) in such a way that one obtains multiple 2-cells, which eventually permit additional arrows between them, giving 3-morphisms as cubes bounded by squares. These cubes then define a 3-fold category as follows:
\begin{definition}
A 3-fold category, denoted   \begin{tikzcd}  D_2 \arrow[r, shift left]   \arrow[r, shift right]  & D_1 \arrow[r, shift left]   \arrow[r, shift right]  & D_0   \end{tikzcd}   is defined in terms of the following conditions:

\begin{enumerate}[(i)]
\item
The objects are the objects of $D_0$ 
\item 
The vertical arrows are the morphisms of $D_0$ 
\item
The horizontal arrows are the objects of $D_1$  
\item
The vertical squares are the morphisms of $D_1$ 
\item
The horizontal squares are the objects of $D_2$  
\item
The cube bounded by vertical and horizontal squares is a morphism of $D_2$ 
\end{enumerate}
The above data has to satisfy the commutative diagram:
\begin{equation}
\begin{tikzcd}[row sep=1.5em]
A \arrow[rr,"f",""{name=0, below}] \arrow[dr,swap,"a",""{name=1}""] \arrow[dd,swap,"h",""{name=2}""] &&
  B \arrow[dd,swap,"h'" near end,""{name=3}""] \arrow[dr,"b",""{name=4, below}""] \\
& A'\arrow[rr, crossing over, "f'" left,""{name=U, below}"", yshift=-0.05ex]
&& B' \arrow[dd,"k'",""{name=5}""] \\
C \arrow[rr,"g" near start,""{name=N, below}""] \arrow[dr,swap,"c",""{name=6}"", below] && D \arrow[dr,"d",""{name=7, below}""] \\
& C' \arrow[rr,"g'" near start,""{name=M}""] \arrow[uu,<-,crossing over,"k" near end]&& D'
 \arrow[from=1,to=4, dash, shorten=10mm, phantom,""{name=K}""]
 \arrow[r,from=6,to=7, dash, shorten=10mm, phantom,""{name=L}""]
 \arrow[Rightarrow, line width=0.3pt, from=K, to=L,"\psi"] \\
 \arrow[r, dash, line width=0.3pt, from=K, to=L]
 \arrow[Rightarrow, from=0, to=N, xshift=2ex, crossing over, shorten=1mm, "\eta", near start, swap]
 \arrow[Rightarrow, from=U, to=M, xshift=2ex, crossing over, shorten=1mm, "\omega", near end, swap]
 \arrow[Rightarrow, from=1, to=6, shorten=2mm, "\alpha",near start]
 \arrow[Rightarrow, from=4, to=7, shorten=1mm, "\beta"]
 \arrow[r, from=1, to=4, shorten =11mm, xshift=10ex, yshift=-3.03ex, crossing over]
 \arrow[r, dash, from=1, to=4, shorten=11mm, xshift=-0.3ex, yshift=-3.03ex, crossing over]
 \end{tikzcd} 
\end{equation}
where arrows are 1-morphisms, squares are 2-morphisms and cubes are 3-morphisms (3-cells) of the 3-fold category.   Compositions, identity and associative laws on each of these entities follow from their usual definitions in $D_0$, $D_1$ and $D_2$ where appropriately applicable.       
\end{definition}
The 3-cells in a multiway system can then be obtained by vertical, horizontal and sideways compositions of cubes of a 3-fold category.

Likewise, one can continue this process of constructing higher morphisms indefinitely, so long as additional rewrite rules between cells of up to order $n - 1$ are admissible such that $n$-morphisms are now $n$-hypercubes defined within an $n$-fold category, and $n$-cells of such multiway systems are simply compositions of $n$-hypercubes in $n$-directions. The ensuing  $n$-fold category can then be expressed as:
\[  \begin{tikzcd} {\cal M}_{n - 1} \arrow[r, shift left] \arrow[r, shift right]  & {\cal M}_{n - 2} \,\,  \cdots\cdots \,\, {\cal M}_2 \arrow[r, shift left] \arrow[r, shift right] & {\cal M}_1 \arrow[r, shift left] \arrow[r, shift right]  & {\cal M}_0 \end{tikzcd} \]
where the objects of ${\cal M}_i$ (for $0 < i < {n-2}$) are $i$-dimensional hypercubes in ${\mathbb R}^n$ whose normal vectors are oriented along the $z$-axis in ${\mathbb R}^n$. The morphisms of ${\cal M}_i$ are $(i + 1)$-dimensional hypercubes in ${\mathbb R}^n$ whose normal vectors are oriented orthogonal to the $z$-axis in ${\mathbb R}^n$. Furthermore, the objects of ${\cal M}_{n - 1}$  are $(n - 1)$-dimensional hypercubes with normal vectors oriented along the $z$-axis in ${\mathbb R}^n$ and the morphism in ${\cal M}_{n - 1}$ is the commutative $n$-hypercube composed of the above. 

This iterative definition of an $n$-fold category corroborates the previously described construction, in which homotopies of up to order $n$ in multiway systems are realized by including additional rewrite rules that introduce new arrows to obtain squares (and so on), until one eventually obtains a commutative $n$-hypercube. 
\qed  \\  \,

The benefit gained by realizing $n$-fold categories in our constructions is the ease of expressing multiple compositions by gluing hypercubes up to order $n$. Composable arrays of hypercubes can then be used to construct all ($i \leq n$)-cells of the category. 

For example, in a 2-fold category, we can define a composable array of 2-dimensional elements (squares) to be such that any array is composable with its immediate neighbors (the associative and commutative laws imply that the composition is well defined). This process easily extends in $n$-fold categories  using $n$-dimensional elements ($n$-hypercubes). 

The following remark is in order:
\begin{remark}
As in Corollary \ref{c3.1}, a $n$-homotopy multiway system extends to an $n$-fold groupoid upon admitting invertible rules that ensure invertibility of all higher morphisms.
\label{rem:Remark1}
\end{remark}

%\begin{remark}
%Even though an $n$-fold category is a strict version of an $n$-category, in that all $n$ composition operations are strictly unital and associative and strictly commute with each other, nonetheless, $n$-fold groupoids model all order-$n$  homotopy types, and this fact is what ensures that multiway rewriting systems are models of order-$n$ homotopy types.
%\end{remark}

%\subsection{The $\infty$-Limit of Multiway Rewriting Systems as $\infty$-Groupoids}

Given the iterative construction above, where one induces higher homotopies in multiway systems via the inclusion of supplementary rules, one can now ask what the $n \to \infty$ limit of this construction yields? 

\begin{proposition}
The $n \to \infty$ limit of a  multiway rewriting system with invertible homotopy rules, corresponding to higher morphisms (so long as they are admissible), is an $\infty$-groupoid. 
\label{p3.3}
\end{proposition}

{\it Proof of Proposition \ref{p3.3}}: 
Following Proposition \ref{p3.2} and Remark \ref{rem:Remark1}, multiway rewriting systems equipped with homotopies up to order $n$, along with invertible rules, yield an $n$-fold groupoid. 
The $n \to \infty$ limit of an $n$-fold groupoid is precisely an $\infty$-groupoid. This limit exists so long as the infinite hierarchy of rules required to iteratively construct higher morphisms is admissible on that multiway rewriting system.  
\qed \\ \, 

%The following remark is due:
%\begin{remark}
%The $\infty$-groupoids discussed here have been constructed within the model of strict n-categories. Even though some definitions of $\infty$-categories (or statements about the homotopy hypothesis based on those definitions) that declare  $\infty$-groupoids as topological spaces are formulated within the model of weak or quasi categories, there is indication that such definitions or statements will eventually be expressible within other models of categories \cite{Riehl2017a}. Hence, it would be interesting to consider what such weakenings would be interpretable as within the kinds of rewriting systems considered above. 
%\end{remark}
%

\section{Implications and Outlook}

This paper has illustrated how higher categorical structures arise naturally within the abstract non-deterministic rewriting systems of the Wolfram model, at least when formalized as multiway systems.  We have shown how higher homotopies induced on multiway systems via specific rewriting rules correspond to morphisms of an $n$-fold category, and we have both demonstrated an explicit computational procedure for constructing these homotopies on multiway rewriting systems, as well as provided the formal correspondence of these systems to $n$-fold categories and $n$-fold groupoids. Interestingly, the $n \to \infty$ limit of this structure yields an $\infty$-groupoid, with the latter being relevant from the point of view of Grothendieck's homotopy hypothesis \cite{baez2007homotopy}; we speculate that this correspondence might potentially offer a new way to understand how the kinds of spatial structures relevant to the foundations of physics might potentially arise from purely combinatorial rewriting systems \cite{Arsiwalla2020}\cite{Shulman2017}\cite{Baez2006}. Our results here provide a preliminary framework for investigating such connections to topological and geometric spaces from an underlying computational perspective, which we outline in forthcoming work \cite{arsiwalla2021pregeometric}.

On the other hand, complementary to what we have shown here, there is an ongoing research program involving the application of methods from homotopy type theory and infinity-category theory, that seeks to formalize, among other things, the mathematical foundations of synthetic spatial structures (in particular, cohesive topological and geometric spaces) \cite{Shulman2017}\cite{Program2013}\cite{Riehl2017a}. Constructions of synthetic geometry are also of great relevance for the foundations of physics. In particular, there are active efforts in the formalization of quantum field theories on a cohesive $\infty$-topos   \cite{Schreiber2012}. In a sense, the Wolfram model can be thought of as a ``constructivist'' realization of some of the above ideas, expressed via a concrete model of computation, namely abstract rewriting systems. The precise connections of Wolfram model multiway rewriting  systems to homotopy types and infinity-categories is also a topic for future work. 

\subsubsection*{Acknowledgments}  The authors would like to thank Stephen Wolfram for his encouragement, and for the many useful conversations and suggestions that helped shape this work. 

%\nocite{*}

\bibliographystyle{eptcs}
\bibliography{ACTHomotopyBibliography.bib}

\end{document}